\documentclass[times]{my-iapress}

\usepackage{moreverb}

\usepackage[colorlinks,bookmarksopen,bookmarksnumbered,citecolor=red,urlcolor=red]{hyperref}

\def\volumeyear{2018}
\def\volumenumber{6}
\def\volumemonth{March}
\usepackage{amsmath,amssymb}
\usepackage{graphicx}

\usepackage{bbm}
\usepackage{subcaption}
\usepackage{ragged2e}

\newcommand{\parc}[1]{\ensuremath{\left(#1\right)}}

\newtheorem{definition}{Definition}

\begin{document}

\runningheads{\ \ \ Parameter Estimation, Sensitivity Analysis
and Optimal Control of a Periodic Epidemic Model}{S. Rosa and D. F. M. Torres}

\title{Parameter Estimation, Sensitivity Analysis and Optimal Control\\
of a Periodic Epidemic Model with Application to HRSV in Florida}

\author{Silv\'erio Rosa\affil{1}, Delfim F. M. Torres\affil{2}$^,$\corrauth}

\address{\affilnum{1}Department of Mathematics and Instituto de Telecomunica\c{c}\~oes (IT),
University of Beira Interior, 6201-001 Covilh\~a, Portugal.\\
\affilnum{2}Center for Research and Development in Mathematics and Applications (CIDMA),
Department of Mathematics, University of Aveiro, 3810-193 Aveiro, Portugal.}

\corraddr{Delfim F. M. Torres (Email: delfim@ua.pt).
Department of Mathematics, University of Aveiro, 3810-193 Aveiro, Portugal.}


\begin{abstract}
A state wide Human Respiratory Syncytial Virus (HRSV) surveillance
system was implemented in Florida in 1999 to support clinical
decision-making for prophylaxis of premature infants.
The research presented in this paper addresses the problem of fitting
real data collected by the Florida HRSV surveillance system
by using a periodic SEIRS mathematical model. A sensitivity
and cost-effectiveness analysis of the model is done
and an optimal control problem is formulated and solved
with treatment as the control variable.
\end{abstract}


\maketitle



\section{Introduction}

Human respiratory syncytial virus (HRSV) is a virus that causes
respiratory tract infections. It is a major cause of lower respiratory
tract infections and hospital visits during infancy and childhood.
A prophylactic medication, palivizumab, can be employed to prevent HRSV in preterm
infants (under 35 weeks gestation), infants with certain congenital heart defects
or bronchopulmonary dysplasia, and infants with congenital malformations
of the airway. Treatment is limited to supportive care, including oxygen therapy.
In temperate climates, there is an annual epidemic during the winter months.
In tropical climates, infection is most common during the rainy season.
In the United States, 60\% of infants are infected during their first HRSV season,
and nearly all children will have been infected with the virus by two to three years of age \cite{Glezen}.
Of those infected with HRSV, 2 to 3\% will develop bronchiolitis, necessitating hospitalization \cite{Hall}.
Natural infection with HRSV induces protective immunity, which wanes over time, possibly
more so than other respiratory viral infections, and thus people can be infected multiple times.
Sometimes an infant can become symptomatically infected more than once, even within a single HRSV season.
Severe HRSV infections have increasingly been found among elderly patients.
Young adults can be re-infected every five to seven years, with symptoms looking
like a sinus infection or a cold. The Florida Department of Health provides an integrated
and reliable HRSV system, with data from hospitals and laboratories \cite{flhealth}.

Mathematical models can project how infectious diseases progress,
to show the likely outcome of an epidemic, and help inform public health
interventions. In epidemiology, compartmental models serve as the base 
mathematical framework for understanding the complex dynamics of these systems. 
Such compartments, in the simplest case, stratify the population into two health states: 
susceptible to the infection of the pathogen, often denoted by $S$, and infected by the 
pathogen, often denoted by the symbol $I$. The way that these compartments interact 
is based upon phenomenological assumptions, and the model is built up from there. 
These models are usually investigated through ordinary differential equations.
To push these basic models to further realism, other compartments are often included, 
most notably the recovered/removed/immune compartment, often denoted by $R$.
A crucial question consists to find parameters
for the particular disease under study, and use those parameters to calculate
the effects of possible control interventions, like treatment or vaccination.
Then the central issue is how to implement such interventions in an optimal way.
This investigation program has been recently carried out for several
infectious diseases, as diverse as dengue \cite{MR3393302,MR3557143},
tuberculosis \cite{MR3388961,MR3562914},
Ebola \cite{MR3578107,MR3544685},
HIV/AIDS \cite{MR3392642,SILVA201770},
and cholera \cite{MR3602689,MR3668112}.
Here we investigate such approach to HRSV.

A comparison of the standard SIRS model
with a more complex model of HRSV transmission,
in which individuals acquire immunity gradually
after repeated exposure to infection, is given in \cite{Weber2001}.
In \cite{MR2718412}, an age-structured mathematical model for HRSV
is proposed, where children younger than one year old,
who are the most affected by this illness, are specially considered.
Real data of hospitalized children in the Spanish region of Valencia
is used in order to determine some seasonal parameters of the model
\cite{MR2718412}. A numerical scheme for the SIRS seasonal epidemiological
model of HRSV transmission is proposed in \cite{MR2435573}. It turns out
that solutions for HRSV compartmental models are typically periodic \cite{MR2426325}.
For this reason, in this work we propose the use of optimal control theory to a
non-autonomous SEIRS model \cite{UBI:UA} and show its usefulness 
in agreement with real HRSV data provided by the Florida Department 
of Health \cite{flhealth}.

The paper is organized as follows. In Section~\ref{sec:2},
we briefly review the SIRS and SEIRS epidemic models.
Our results are then given in Section~\ref{sec:3}:
parameter estimation of the SEIRS model with real data of Florida
(Section~\ref{sec:fit}); sensitivity analysis (Section~\ref{sensitivity});
and optimal control and numerical simulations
(Sections~\ref{section:OC} and \ref{sec:NumericalRes}).
We end with Section~\ref{sec:conc} of conclusions
and some future perspectives.


\section{Review of some recent HRSV mathematical models}
\label{sec:2}

We focus on compartmental models that divide the population into mutually
exclusive distinct groups (of susceptible, or infected, or immune individuals,
or \ldots) and we use  deterministic continuous transitions between those groups,
also known as states. Due to the seasonality of HRSV, the models that best fit
real data are periodic and, therefore, non-autonomous. In \cite{Weber2001},
two models are proposed where the transmission is periodic:
(i) a simple model with only three compartments, known as SIRS (Section~\ref{subsec:SIRS});
(ii) and a more complex model with seventeen compartments, named MSEIRS4. However,
it is shown that the simpler model is closer to real data \cite{Weber2001}.
Zang et al. \cite{zhang2012existence} use a non-autonomous SEIR model where,
beyond the periodicity in the transmission rate,  the annual recruitment
rate is also periodic. This assumption is due to opening and closing
of schools \cite{zhang2012existence}. Here we adopt such ideas and consider
a simple non-autonomous and periodic SEIRS model (Section~\ref{subsec:SEIRS}).


\subsection{SIRS model}
\label{subsec:SIRS}

In the SIRS model, one considers that the population consists
of susceptible ($S$), infected and infectious ($I$),
and recovered ($R$) individuals. A characteristic feature
of HRSV is that immunity after infection is temporary,
so that the recovered individuals become susceptible again \cite{Weber2001}.
The model depends on several parameters: parameter $\mu$, which denotes
the birth rate, assumed equal to the mortality rate;
$\gamma$, which is the rate of immunity loss;
and $\nu$, which is the rate of loss of infectiousness.
Moreover, the influence of the seasonality on the transmission
parameter $\beta$ is modeled by the \emph{cosine} function.
Precisely, and using the linear mass action law \cite{Weber2001}, 
the SIRS model for HRSV is given by the following system of ordinary
differential equations \cite{Weber2001}:
\begin{equation}
\label{eq:modSIRS}
\begin{cases}
\dot{S}(t)=  \mu-\mu S(t)-\beta(t)S(t) I(t) +\gamma R(t),\\
\dot{I}(t)= \beta(t)S(t) I(t) -\nu I(t)-\mu I(t),\\
\dot{R}(t)=\nu I(t)-\mu R(t)-\gamma R(t),
\end{cases}
\end{equation}
where  $\beta(t)=b_0(1+b_1\cos(2 \pi t +\Phi))$.
Note that $b_0$ is the average of the transmission parameter
$\beta$, $b_1$ is the amplitude of the seasonal
fluctuation in the transmission parameter $\beta$,
while $\Phi$ is an angle that will be chosen later in agreement 
with the real data for $I(t)$.


\subsection{SEIRS model}
\label{subsec:SEIRS}

In order to incorporate some important features of HRSV into the model
while keeping it simple, we extend the model \eqref{eq:modSIRS}
of Section~\ref{subsec:SIRS} by using the ideas in \cite{UBI:UA}:
firstly, we include a latency period by introducing a group $E$
of individuals who have been infected but are not yet infectious,
these individuals becoming infectious at a rate $\varepsilon$
and the latency period being equal to the time between
infection and the first symptoms; secondly, and similar
to \cite{zhang2012existence}, we consider that the annual
recruitment rate is seasonal due to schools opening/closing periods.
The equations of the SEIRS model are then given by 
\begin{equation}
\label{eq:modSEIRS}
\begin{cases}
\dot{S}(t) =  \lambda(t)-\mu S(t)-\beta(t)S(t) I(t) +\gamma R(t),\\
\dot{E}(t) = \beta(t)S(t) I(t) -\mu E(t)-\varepsilon E(t),\\
\dot{I}(t) =\varepsilon E(t)-\mu I(t)-\nu I(t),\\
\dot{R}(t) =\nu I(t)-\mu R(t)-\gamma R(t),
\end{cases}
\end{equation}
where $\lambda(t)=\mu(1 + c_1  \cos( 2 \pi t +\Phi) )$
is the recruitment rate, including newborns, immigrants, etc.
Parameter $c_1$ is the amplitude of the seasonal fluctuation
in the recruitment  parameter $\lambda$ and $\Phi$ an angle
to be conveniently chosen.


\section{Main results}
\label{sec:3}

We begin by investigating how realistic the models discussed in
Section~\ref{sec:2} are with respect to HRSV and real data from
Florida \cite{flhealth}. For that we need a proper estimation
of the parameter values.


\subsection{Parameter estimation}
\label{sec:fit}

\begin{figure}[!htb]
\centering
\includegraphics[scale=0.7]{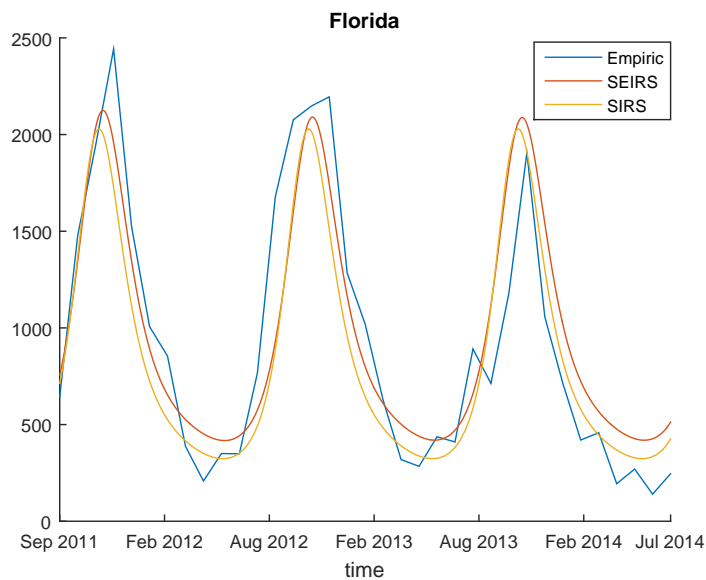}
\caption{Comparison of the number of infective individuals
registered in the state of  Florida \cite{flhealth}
with the ones predicted by models SIRS \eqref{eq:modSIRS}
and SEIRS \eqref{eq:modSEIRS} with the parameter values of
Table~\ref{tab:param}.}
\label{fig:clients}
\end{figure}

During model fitting, the values of the parameters $\mu$ (birth rate), 
$\varepsilon$ (infectious rate), $\nu$ (rate of infectiousness loss) 
and $\gamma$ (rate of immunity loss) were held constant. See Table~\ref{tab:param},
where we use $\overline{n}$ to denote the average number of monthly HRSV cases reported.
The values of $\varepsilon$, $\nu$ and $\gamma$ were taken from \cite{Weber2001}.
For the birth rate, $\mu$, we used the value given in \cite{flhealthcharts}
for the state of Florida. In agreement with Section~\ref{subsec:SIRS},
we have set the birth rate equal to the mortality rate, so that one obtains
a constant population during the period under study. This simplification is justified when,
as in this case, the annual infection rate is much bigger than the growth of population.
To simplify analysis, we set the total population equal to one and use
fractions for the values of state variables. Values of four parameters were determined
by fitting the model: (i) the mean of the transmission parameter, $b_0$;
(ii) its relative seasonal amplitude, $b_1$; (iii) the angle $\Phi$;
and (iv) a scaling factor $s$, which scales the number of infectious individuals
to the empirical case reported in our unit scaled model. The angle $\Phi$ was normalized
in the following way: for a value of zero, a maximum of the cosine function used
in $\beta(t)$ (and also used  in $\lambda(t)$) coincides with
the first maximum of the empirical cases.
The models were fitted to data on the reported number of positive tests of HRSV disease,
per month, in the State of Florida, excluding North region, between September 2011
and July 2014 (35 months). The data was obtained
from the Florida Department of Health \cite{flhealth}.
Our results are given in Figure~\ref{fig:clients}.
\begin{table}[!htb]\caption{Results of model fitting.}
\label{tab:param}
\centering
\begin{tabular}{lccccccccccc}\toprule
model & $\mu$ & $\nu$ & $\gamma$ & $\varepsilon$ & $b_0$ 
& $b_1$ & $c_1$ & $s$ & $\Phi$ & $\overline{n}$ & $R_0$\\[1mm] \midrule
SIRS & 0.0113 & 36 & 1.8 & --  & 74.2 & 0.14 & --  & 35000 & $7\pi/5$ & 932 & 2.06\\[1mm]
SEIRS & 0.0113 & 36 & 1.8 & 91  & 88.25 & 0.17 & 0.17 & 35000 & $7\pi/5$ & 932 & 2.45\\[1mm] \bottomrule
\end{tabular}
\end{table}
Following \cite{silva2017modeling}, our fitting approach consisted to minimize the $l_2$ norm
of the difference between the real data and predictive cases of HRSV infection
given by models \eqref{eq:modSIRS} and \eqref{eq:modSEIRS}.
Let $e$ represent the relative error, computed by
$$
e=\frac{\|I_{\text{model}}-I_{\text{empiric}}\|_2}{\|I_{\text{empiric}}\|_2},
$$
where $I_{\text{empiric}}$ is the real data obtained from \cite{flhealth}
and $I_{\text{model}}$ the one predicted by the SIRS or SEIRS models.
The results shown in Figure~\ref{fig:clients} correspond to a relative error
of 0.066\% and 0.060\% of infants per year with respect to the total
child population of Florida in 2014, excluding North region, respectively
for SIRS or SEIRS models. We note that the values of $\mu$, $\nu$, $\gamma$
and $\varepsilon$ are annual rates. The difference between the two models,
in terms of absolute errors, clearly justifies SEIRS' superiority.


\subsection{Sensitivity analysis}
\label{sensitivity}

One of the most important thresholds while studying
infectious disease models is the basic reproduction
number \cite{VANDENDRIESSCHE2017}. The basic reproduction number
of the SIRS model was computed in \cite{Weber2001}.
Assuming that the transmission and the recruitment parameters are constant,
the basic reproduction number for the SEIRS model is given by (see \cite{MR2559058})
\begin{equation}
\label{r0:SEIRS}
R_0=\frac{\beta\varepsilon}{(\mu+\nu)(\varepsilon+\mu)}.
\end{equation}

Now we do a sensitivity analysis
for the basic reproduction number \eqref{r0:SEIRS}.
Such analysis tells us how important each parameter is to disease transmission.
This information is crucial not only for experimental design,
but also to data assimilation and reduction of complex nonlinear models
\cite{powell2005sensitivity}. Sensitivity analysis is commonly used
to determine the robustness of model predictions to parameter values,
since there are usually errors in data collection and presumed parameter values.
It is used to discover parameters that have a high impact on $R_0$
and should be targeted by intervention strategies. More precisely,
sensitivity indices's allows to measure the relative change
in a variable when parameter changes. For that we use the normalized
forward sensitivity index of a variable, with respect to a given parameter,
which is defined as the ratio of the relative change in the variable
to the relative change in the parameter. If such variable is differentiable
with respect to the parameter, then the sensitivity index
is defined using partial derivatives, as follows
(see \cite{chitnis2008determining,rodrigues2013sensitivity}).

\begin{definition}
\label{def:sentInd}
The normalized forward sensitivity index of $R_0$, which is differentiable
with respect to a given parameter $p$, is defined by
$$
\Upsilon_p^{R_0}=\frac{\partial R_0}{\partial p}\frac{p}{R_0}.
$$
\end{definition}

\begin{figure}[!htb]
\centering
\begin{subfigure}[b]{0.46\textwidth}\centering
\includegraphics[scale=0.46]{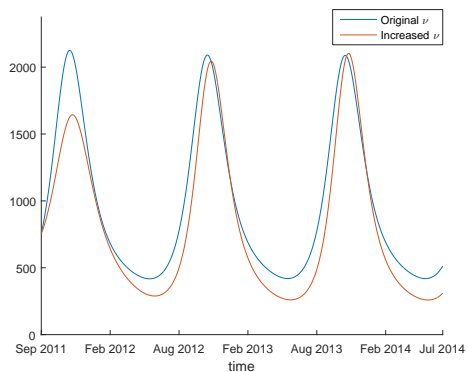}
\caption{Impact of the variation of $\nu$ in the number of infectives $I$.}
\label{fig:sensitivity_nu}
\end{subfigure}
\hspace*{1cm}
\begin{subfigure}[b]{0.46\textwidth}\centering
\includegraphics[scale=0.46]{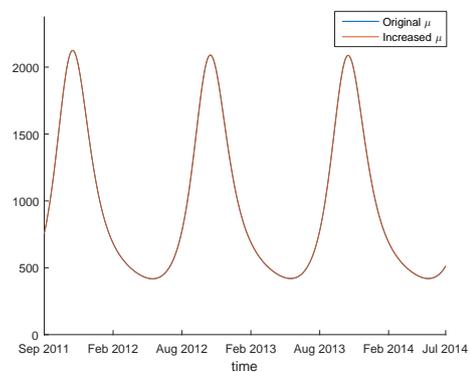}
\caption{Impact of the variation of $\mu$ in the number of infectives $I$
(difference not visible).}
\label{fig:sensitivity_miu}
\end{subfigure}
\caption{Infected number of individuals predicted by SEIRS model
with original parameter values as in Table~\ref{tab:param}
and with an increase of 10\% of a specific parameter:
$\nu$ in (a) and $\mu$ in (b).}
\label{fig:sensitivity}
\end{figure}

One can easily compute an analytical expression for the sensitivity of $R_0$,
using the explicit formula \eqref{r0:SEIRS}, to each parameter that it includes.
The values of the sensitivity indices for the parameters values of Table~\ref{tab:param}
are presented in Table~\ref{tab:sensitivity}. Note that the sensitivity index
may depend on several parameters of the system, but also can be constant,
independent of any parameter. For example, $\Upsilon_{\beta}^{R_0}=+1$,
meaning that increasing (decreasing) $\beta$ by a given percentage increases
(decreases) always $R_0$ by that same percentage.
\begin{table}[!htb]
\centering
\caption{Sensitivity of $R_0$ evaluated for the parameter values given
in Table~\ref{tab:param}, according with \eqref{r0:SEIRS} and
Definition~\ref{def:sentInd}.}\label{tab:sensitivity}
\begin{tabular}{cc@{\hspace*{5cm}}c}\toprule
Parameter && Sensitivity index \\[1mm] \midrule
$\beta$ && $+1$\\[1mm]
$\varepsilon$ && $+0.283465$\\[1mm]
$\nu$ && $-1.28315$\\[1mm]
$\mu$ && $-0.00031379$\\ \bottomrule
\end{tabular}
\end{table}
The estimation of a sensitive parameter should be carefully done,
since a small perturbation in such parameter leads to relevant
quantitative changes. On the other hand, the estimation of
a parameter with a small value for the sensitivity index
does not require as much attention to estimate,
because a small perturbation in that parameter leads
to small changes \cite{mikucki2012sensitivity}.
According with Table~\ref{tab:sensitivity}, we should pay
special attention to the estimation of parameters $\nu$, $b_0$
(mean value of $\beta(t)$) and $\varepsilon$. In contrast,
the estimation of the rate of birth, $\mu$, (or rate of mortality)
does not require as much attention because of its low value of the sensitivity index.
This is well illustrated in Figure~\ref{fig:sensitivity},
where we can see in Figure~\ref{fig:sensitivity_nu}
the graphics of the number of infectives with and without
an increment of 10\% for the parameter $\nu$ (most sensitive parameter),
and in Figure~\ref{fig:sensitivity_miu} for parameter $\mu$
(less sensitive parameter).


\subsection{Optimal control}
\label{section:OC}

To investigate some optimal control measures,
we choose the SEIRS model, which provides
a better fitting of available real data,
as shown in Section~\ref{sec:fit}. The evolution
on the number of susceptible, exposed, infective and recovered,
depend on some factors that can be controlled. In case of
HRSV, treatment is the most realistic. For this reason,
we consider the following optimal control problem:
\begin{equation}
\label{eq:modSEIRS_control}
\begin{cases}
\dot{S}(t) = \lambda(t)-\mu S(t)-\beta(t)S(t) I(t) +\gamma R(t),\\
\dot{E}(t) = \beta(t)S(t) I(t) -\mu E(t)-\varepsilon E(t),\\
\dot{I}(t) = \varepsilon E(t)-\mu I(t)-\nu I(t)-\mathbbm{T}(t)I(t),\\
\dot{R}(t) = \nu I(t)-\mu R(t)-\gamma R(t)+\mathbbm{T}(t)I(t),
\end{cases}
\end{equation}
subject to given initial conditions
\begin{equation}
\label{eq:cond_ini}
S(0), E(0), I(0), R(0)\geqslant 0,
\end{equation}
where $\mathbbm{T}$ denotes \emph{treatment} ($\mathbbm{T}$ is the control variable).
Note that in the absence of control measures, that is, when
$\mathbbm{T}(t) \equiv 0$, system \eqref{eq:modSEIRS_control}
reduces to \eqref{eq:modSEIRS}.
We consider the set of admissible control functions
\begin{equation}\notag
\Omega=\left\{\mathbbm{T}(\cdot)\in L^{\infty}(0,t_f):
0\leqslant \mathbbm{T}(t)\leqslant \mathbbm{T}_{\text{max}},
\forall t\in[0,t_f]\right\}.
\end{equation}
Our purpose is to minimize the number of infectious individuals
and the cost required to control the disease by treating the infected.
Precisely, the optimal control problem consists in
\begin{equation}
\label{cost-functional}
\begin{split}
\min ~\mathcal{J}(I,\mathbbm{T})
&=\int_0^{t_f} \left(\kappa_1\,I(t)+ \kappa_2\, \mathbbm{T}^2(t)\right) ~dt\\
\mbox{subject to }
&\dot{S}(t) =  \lambda(t)-\mu S(t)-\beta(t)S(t) I(t) +\gamma R(t)\\
&\dot{E}(t) = \beta(t)S(t) I(t) -\mu E(t)-\varepsilon E(t)\\
&\dot{I}(t) =\varepsilon E(t)-\mu I(t)-\nu I(t)-\mathbbm{T}(t) I(t)\\
&\dot{R}(t) =\nu I(t)-\mu R(t)-\gamma R(t)+\mathbbm{T}(t) I(t)
\end{split}
\end{equation}
with $0<\kappa_1,\kappa_2 <\infty$. To solve the problem,
we apply the celebrated Pontryagin maximum principle \cite{MR898009}:
the Hamiltonian is given by
\begin{multline*}
\mathcal{H}
= \kappa_1 I +\kappa_2 \mathbbm{T}^2+p_1(\lambda-\mu S-\beta S I +\gamma R)
+p_2(\beta S I -\mu E-\varepsilon E)\\
+p_3(\varepsilon E-\mu I-\nu I-\mathbbm{T}I)
+p_4(\nu I-\mu R-\gamma R+\mathbbm{T}I);
\end{multline*}
the maximality condition ensures that the extremal control is given by
\begin{equation}
\label{eq:ext:cont}
\mathbbm{T}(t)=\min\left\{\max\left\{0,\dfrac{(p_3(t)-p_4(t))
I(t)}{2 \kappa_2}\right\},\mathbbm{T}_{\max}\right\};
\end{equation}
while the adjoint system asserts that the co-state variables
$p_i\parc{t}$, $1\leq i\leq 4$, satisfy relations
\begin{align}\begin{split}
\dot p_1  =-\frac{\partial \mathcal{H}}{\partial S} \Leftrightarrow \dot p_1
&=p_1(\mu+\beta(t)I)-\beta(t) I p_2,\\[2mm]
\dot p_2  =-\frac{\partial \mathcal{H}}{\partial E} \Leftrightarrow \dot p_2
&= p_2(\mu+\varepsilon)-\varepsilon p_3, \\[2mm]
\dot p_3  = -\frac{\partial \mathcal{H}}{\partial I} \Leftrightarrow \dot p_3
&= -\kappa_1+\beta(t) p_1 S-p_2 \beta(t) S+p_3(\mu+\nu+\mathbbm{T})-p_4(\nu+\mathbbm{T}),\\[2mm]
\dot p_4 =  -\frac{\partial \mathcal{H}}{\partial R} \Leftrightarrow \dot p_4
&= -\gamma p_1+p_4(\mu+\gamma)\end{split}\label{eq:co_states}
\end{align}
and transversality conditions
\begin{equation}
\label{eq:TC}
p_1\parc{t_f}= p_2\parc{t_f}= p_3\parc{t_f}= p_4\parc{t_f}= 0.
\end{equation}
Therefore, in order to solve the optimal control problem
\eqref{cost-functional}, we solve the following boundary value problem:
\begin{equation}
\left\{\begin{aligned}
\dot S =&  \lambda(t)-\mu S-\beta(t)S I +\gamma R\\[1mm]
\dot E =& \beta(t)S I -\mu E-\varepsilon E\\[1mm]
\dot I =&\varepsilon E-\mu I-\nu I-\mathbbm{T}I\\[1mm]
\dot R =&\nu I-\mu R-\gamma R+\mathbbm{T}I\\[1mm]
\dot p_1 =&p_1(\mu+\beta(t)I)-\beta(t) I p_2 \\[2mm]
\dot p_2 =& p_2(\mu+\varepsilon)-\varepsilon p_3 \\[2mm]
\dot p_3 =& -\kappa_1+\beta(t) p_1 S-p_2 \beta(t) S+p_3(\mu+\nu+\mathbbm{T})-p_4(\nu+\mathbbm{T}) \\[2mm]
\dot p_4 =& -\gamma p_1+p_4(\mu+\gamma)
\end{aligned}\right. \label{eq:bvp}
\end{equation}
with given initial conditions
$$
S(0) =s_0,\quad E(0) =e_0,\quad I(0) =i_0,\quad R(0) =r_0
$$
and transversality conditions \eqref{eq:TC},
where $\mathbbm{T}(t)$ is given by \eqref{eq:ext:cont}.
This is done numerically in Section~\ref{sec:NumericalRes}.


\subsection{Numerical results and cost-effectiveness analysis for the optimal control model}
\label{sec:NumericalRes}

Two algorithms were implemented to obtain and confirm the numerical results.
One approach solves the optimal control problem numerically using a fourth
order Runge--Kutta iterative method. First we solve the system
\eqref{eq:modSEIRS_control} with initial conditions for the state variables
and a guess for the control over the time interval $[0,t_f]$, by the forward
Runge--Kutta fourth order procedure, and obtain the values of the state  variables $S$,
$E$, $I$ and $R$. Using those values, then we solve the system \eqref{eq:co_states}
with the transversality conditions \eqref{eq:TC}, by the backward fourth order
Runge--Kutta procedure, and obtain the values of the co-state variables. The
control is updated by a convex combination of the previous control and the value
from \eqref{eq:ext:cont}. The iteration is stopped when the values of the
unknowns at the earlier iteration are very close to the ones at the current
iteration. The results coincide with the ones obtained by a similar iterative
method that, in each iteration, solves the boundary value problem \eqref{eq:bvp}
using the \textsc{Matlab} \textit{bvp4c} routine.

In what follows, we consider that, for simplicity, $\mathbbm{T}_{\max}=1$ 
and the other parameters are fixed according to Table~\ref{tab:param},
with exception to the angle $\Phi$ that is assumed to be $\pi/2$.
Such value allows that the transmission parameter initial value be the
average, $\beta(0)=b_0$, and the recruitment rate initial value 
also be the average, $\lambda(0)=\mu$. The initial conditions, 
given by Table~\ref{tab:solinit}, are obtained as the
nontrivial equilibria for the system \eqref{eq:modSEIRS_control} with no
control, corresponding to the population state prior the introduction of the
treatment. We also assume that $t_f=5$ because the  World Health Organization
goals for most diseases are usually fixed for 5 years periods.
\begin{table}[!htb]
\centering
\caption{Initial conditions for the optimal control problem
with parameters given by Table~\ref{tab:param}, excepting
angle $\Phi$ that is assumed to be $\pi/2$.
The values correspond to the endemic equilibrium of \eqref{eq:modSEIRS_control}
before the introduction of treatment.}\label{tab:solinit}
\begin{tabular}{c@{\hspace*{2.2cm}}c@{\hspace*{2.2cm}}c@{\hspace*{2.2cm}}c}\toprule
$s S(0)$ & $s E(0)$ & $s I(0)$ & $s R(0)$ \\[1mm] \midrule
$14284$  & $385$ &  $974$ &  $19357$ \\ \bottomrule
\end{tabular}
\end{table}
\begin{figure}[!htb]
\centering
\begin{subfigure}[b]{0.46\textwidth}\centering
\includegraphics[scale=0.46]{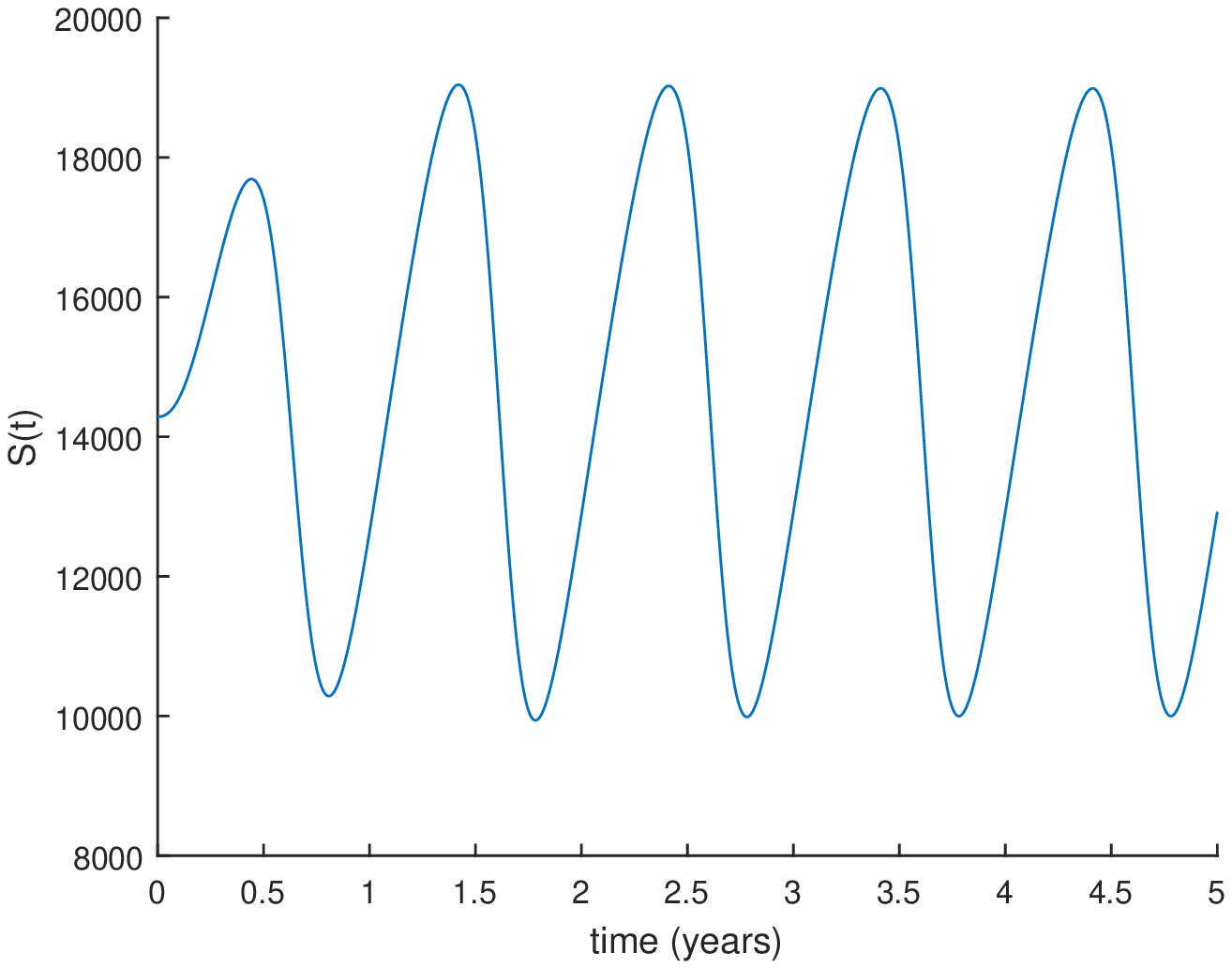}
\caption{Evolution of the number of susceptible individuals.}
\end{subfigure}
\hspace*{1cm}
\begin{subfigure}[b]{0.46\textwidth}
\centering
\includegraphics[scale=0.46]{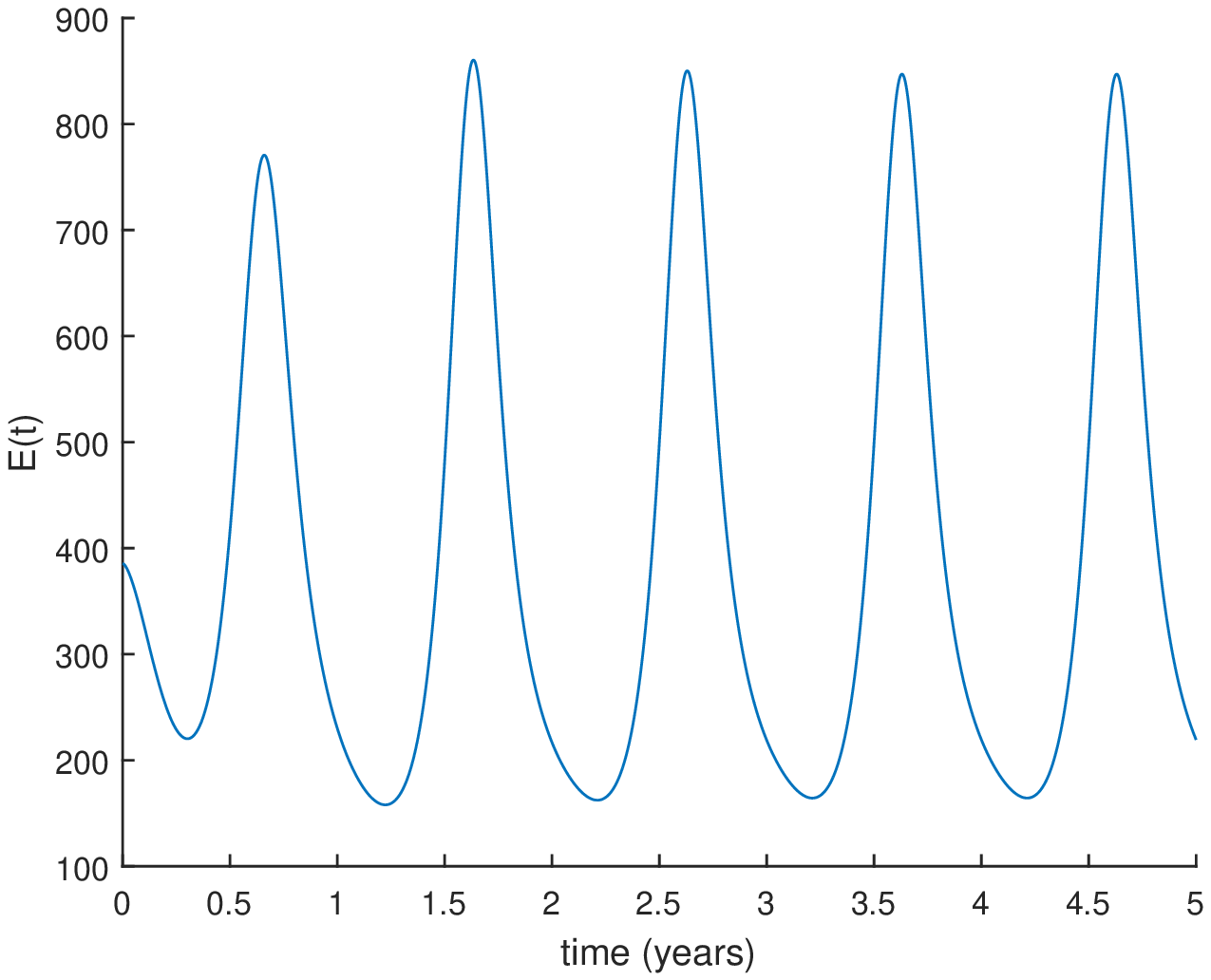}
\caption{Evolution of the number of exposed individuals.}
\end{subfigure}\\
\begin{subfigure}[b]{0.46\textwidth}
\centering
\includegraphics[scale=0.46]{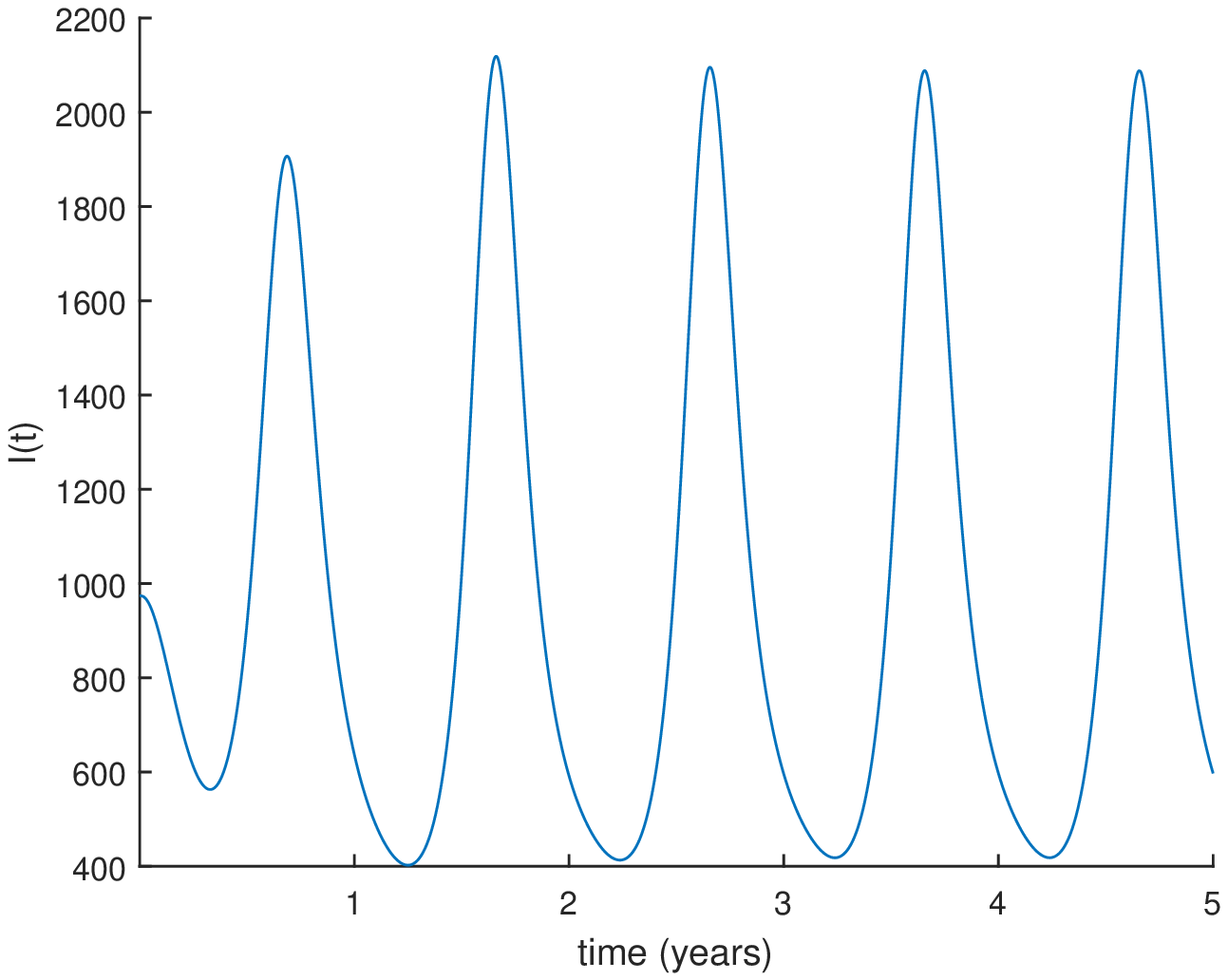}
\caption{Evolution of the number of infected individuals.}
\end{subfigure}\hspace*{1cm}
\begin{subfigure}[b]{0.46\textwidth}
\centering
\includegraphics[scale=0.46]{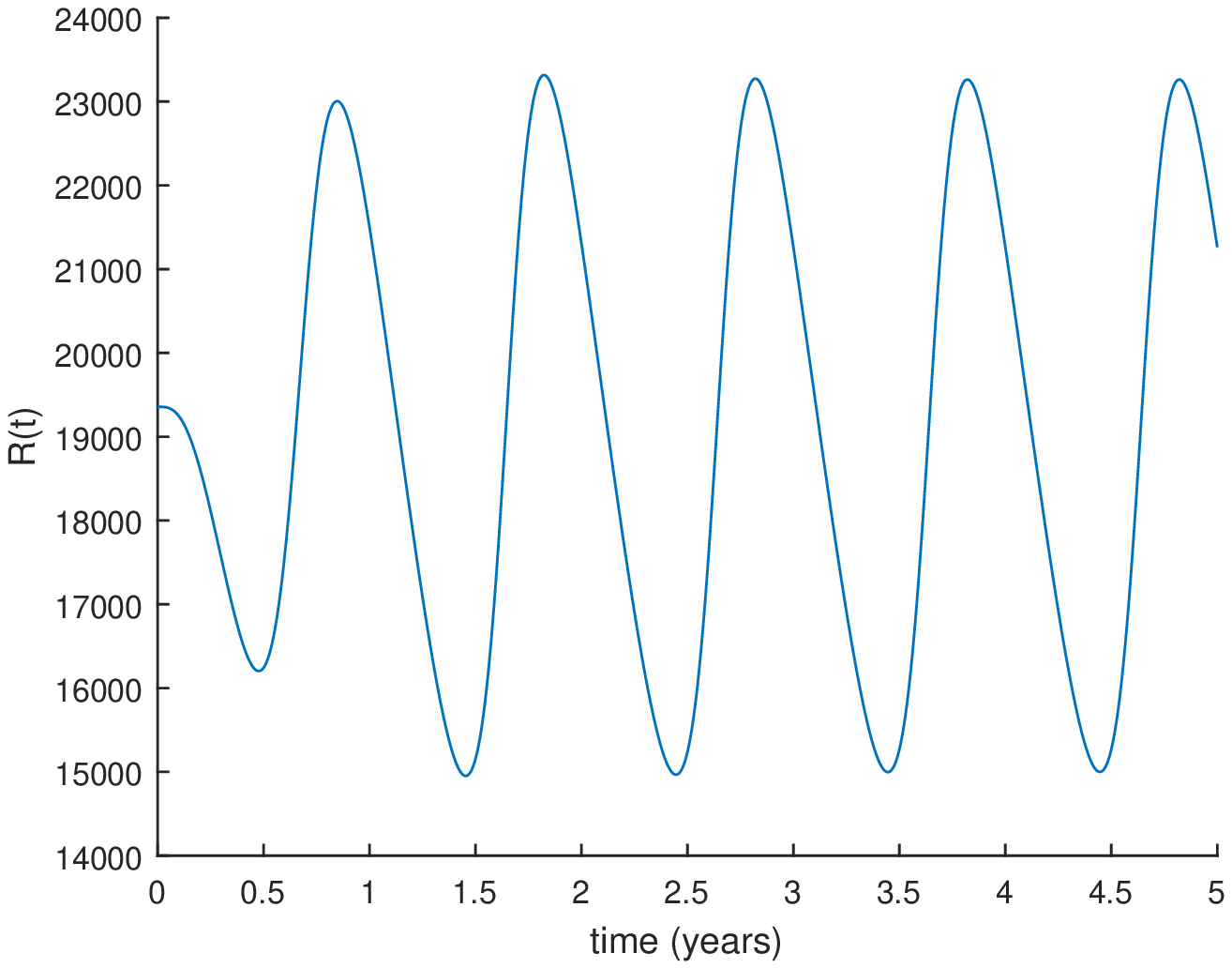}
\caption{Evolution of the number of recovered individuals.}
\end{subfigure}
\caption{State variables of the optimal control problem \eqref{cost-functional},
assuming weights $k_1=1$ and $k_2=0.001$.}
\label{fig:states_var:K}
\end{figure}
\begin{figure}[!htb]
\centering
\begin{subfigure}[b]{0.46\textwidth}\centering
\includegraphics[scale=0.46]{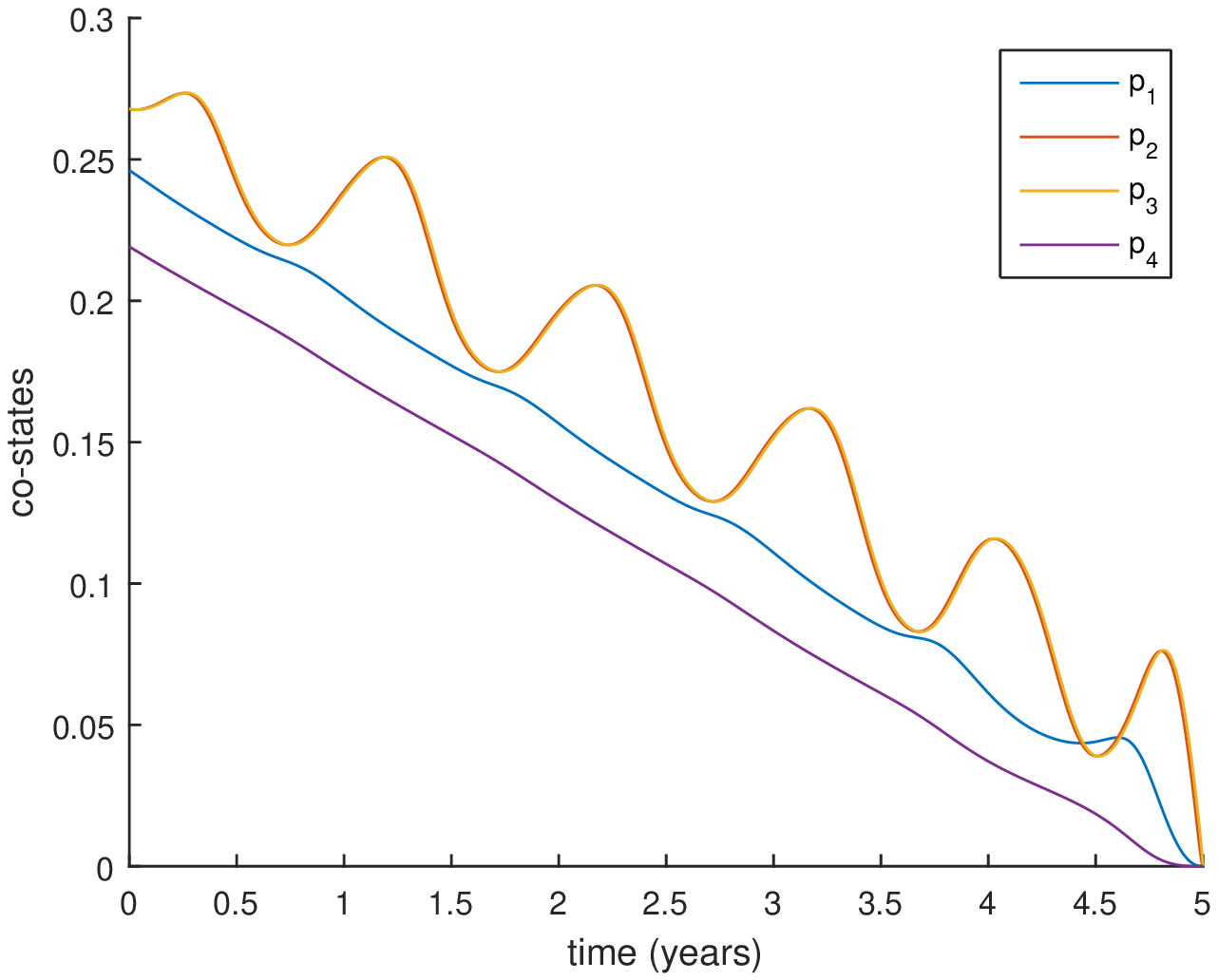}
\caption{Variation of co-state variables.}
\label{fig:co_states_001}
\end{subfigure}
\hspace*{1cm}
\begin{subfigure}[b]{0.46\textwidth}
\centering
\includegraphics[scale=0.46]{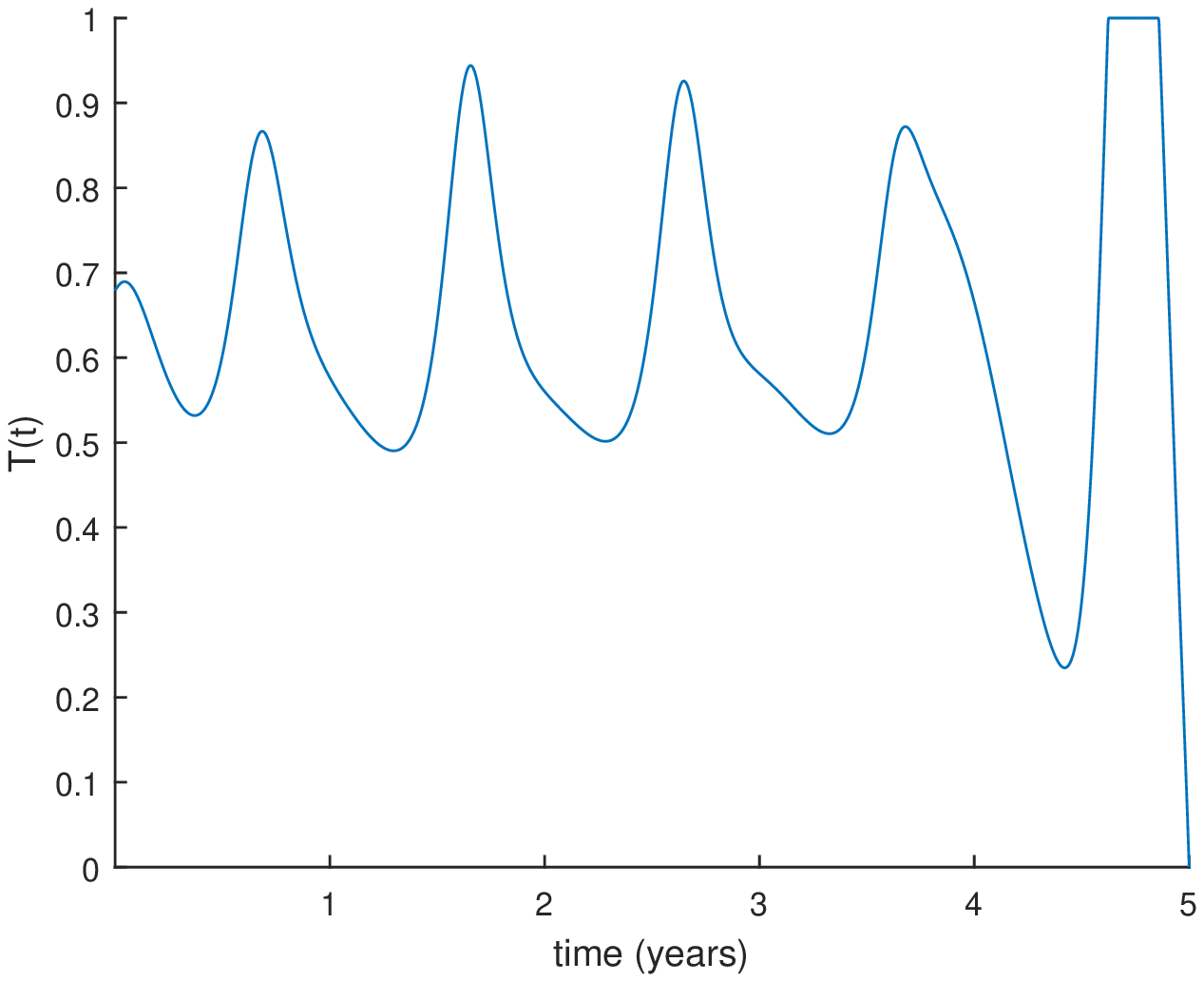}
\caption{Variation of the optimal control $\mathbbm{T}$ (treatment).}
\label{fig:control_001}
\end{subfigure}\\
\begin{subfigure}[b]{0.46\textwidth}
\centering
\includegraphics[scale=0.46]{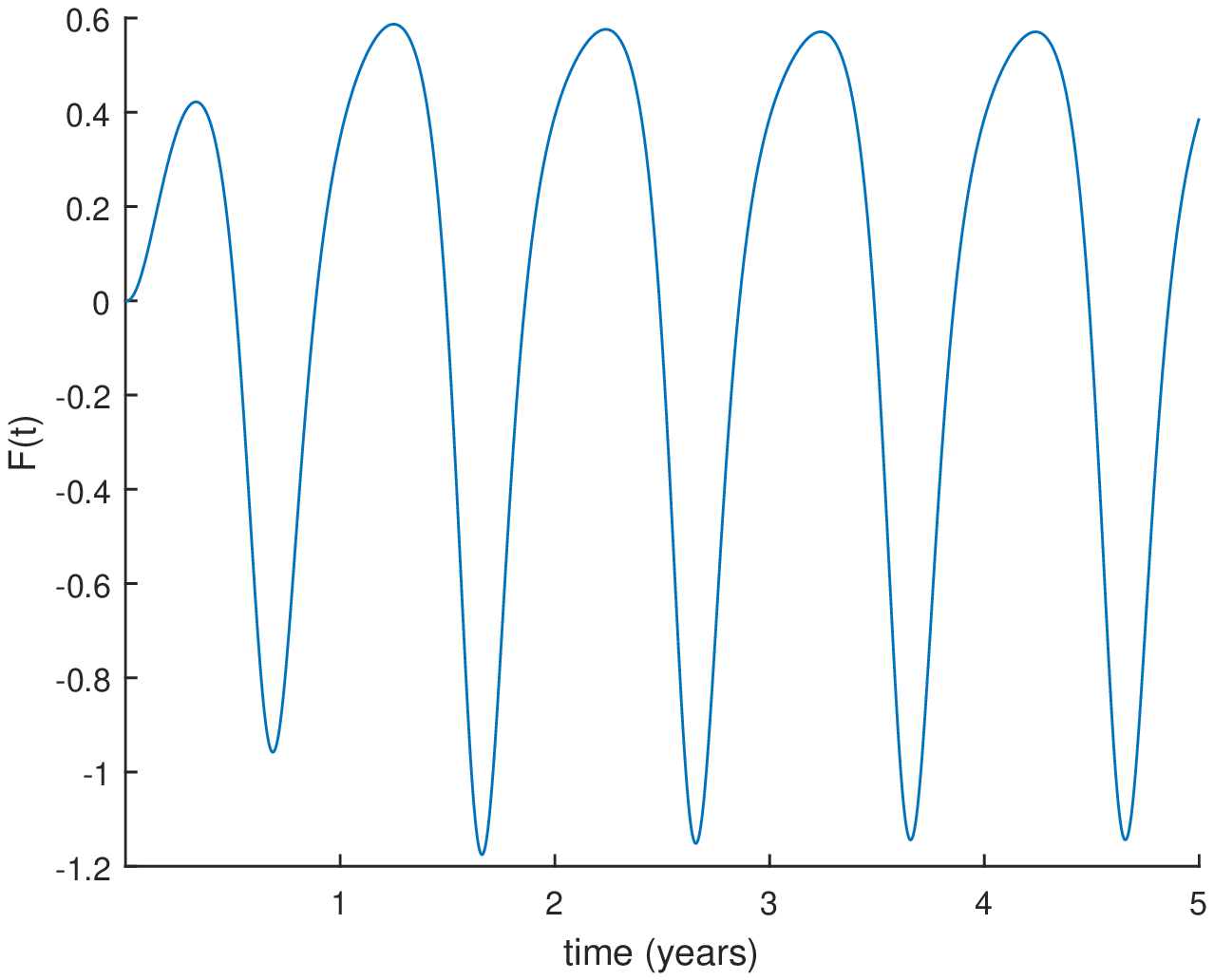}
\caption{Efficacy function $F(t)$ defined in \eqref{efficacy_function}.}
\label{fig:efficacy_001}
\end{subfigure}
\caption{Co-state variables, optimal control $\mathbbm{T}$ and efficacy
function $F(t)$ associated to the optimal control problem
\eqref{cost-functional},
assuming weights $k_1=1$ and $k_2=0.001$.}
\label{fig:control_efficacy:K}
\end{figure}

The solution for the optimal control problem is illustrated in 
Figures~\ref{fig:states_var:K}, \ref{fig:co_states_001} and \ref{fig:control_001}.
The periodic nature of the disease influences the evolution of the four
state variables. We can also see that the control is a continuous function
that shows some nonregularity at the end of the time interval $[0,t_f]$. This
behaviour is explained by the irregular oscillation of the co-state
variables, on which the control depends.

Treatment intensity of the infectious individuals must have, in each
year of the time interval, a given period of time, during which most 
of the infectious individuals are treated, to ensure that the level of
infectious reach very low levels. Figure~\ref{fig:efficacy_001}
shows the efficacy function \cite{rodrigues2014cost}
that is defined by
\begin{equation}
\label{efficacy_function}
F(t)=\frac{I(0)-I^*(t)}{I(0)}=1-\frac{I^*(t)}{I(0)},
\end{equation}
where $I^*(t)$ is the optimal solution associated with the optimal
control and $I(0)$ is the corresponding initial condition. This
function measures the proportional variation  in the number of
infectious individuals after the intervention with control
$\mathbbm{T}^*$, by comparing the number of infected individuals
at time $t$ with the initial value $I(0)$ for which there is no
control implemented. We observe that $F(t)$ oscillates between
$-1.18$ (lower bound) and $+0.59$ (upper bound), and  exhibits
the inverse tendency of $I(t)$.

\begin{figure}[!htb]
\centering
\begin{subfigure}[b]{0.46\textwidth}
\includegraphics[scale=0.46]{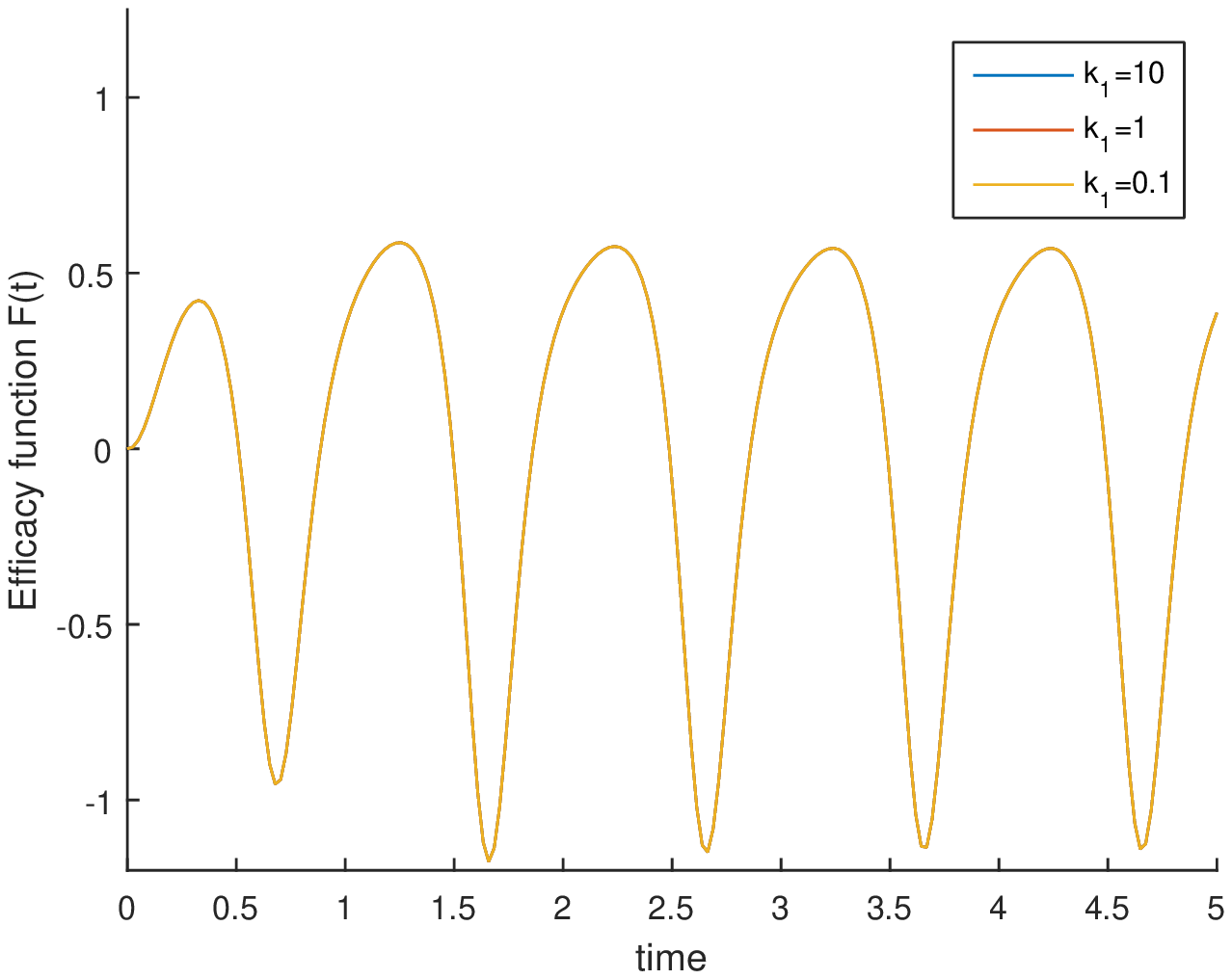}
\caption{$k_2=0.001$ and varying $k_1$.}
\label{fig:infeted_001b}
\end{subfigure}
\hspace*{1cm}
\begin{subfigure}[b]{0.46\textwidth}
\includegraphics[scale=0.46]{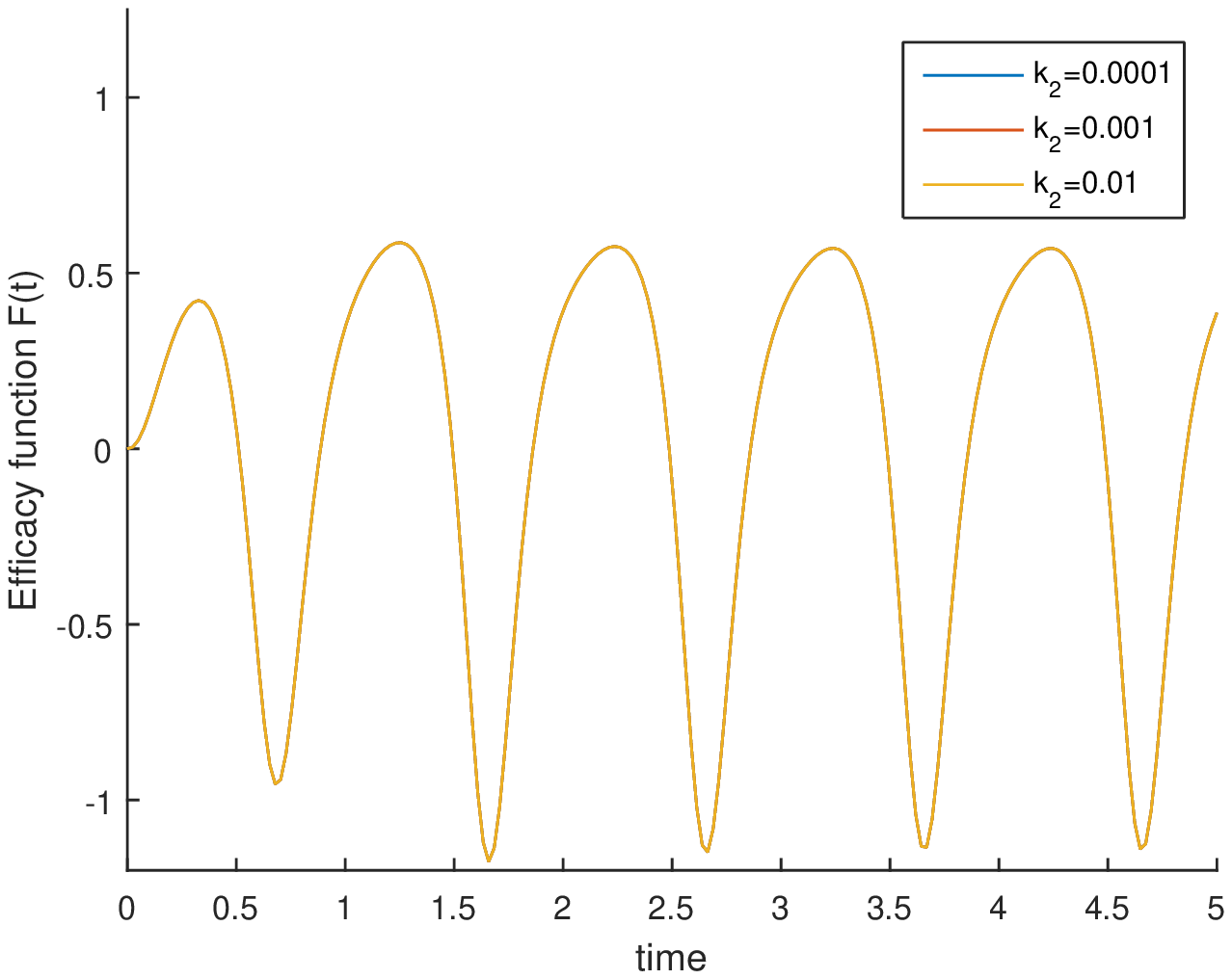}
\caption{$k_1=1$ and varying $k_2$.}
\label{fig:efficacy_001b}
\end{subfigure}
\caption{Sensitivity analysis for the weights of the objective
functional on \eqref{cost-functional}. \emph{Left:} $k_2=0.001$
and $k_1=0.1,1,10$ \emph{Right:} $k_1=1$ and $k_2=0.01,0.001,0.0001$.}
\label{fig:sensitivity:OF}
\end{figure}

Obviously, the results depend on the objective functional
$\mathcal{J}$ given in \eqref{cost-functional}. Namely, they
depend on the weight constants associated with the amount of
infectious individuals $k_1$ and with the cost of the control
$k_2$. According with Figure~\ref{fig:sensitivity:OF}, results
do not change qualitatively by varying constants $k_i, i=1,2$.
Nevertheless, the magnitude of the efficacy changes
slightly when $k_1$ and $k_2$ vary independently.

Some summary measures are introduced to evaluate the cost and
the effectiveness of the proposed control measure for the
intervention period.

The total cases averted by the intervention during
the time period $t_f$ \cite{rodrigues2014cost} is given by
\begin{equation}
\label{eq:A}
A=t_f I(0)-\int_0^{t_f}I^*(t)~dt,
\end{equation}
where $I^*(t)$ is the optimal solution associated with the optimal control
$\mathbbm{T}^*$ and $I(0)$ is the corresponding initial condition. Note
that this initial condition is obtained as the equilibrium proportion
$\overline{I}$ of system \eqref{eq:modSEIRS_control} without treatment
intervention, which is independent on time, so 
$t_f I(0)=\int_0^{t_f}\overline{I}~dt$ represents the total infectious
cases over a given period of $t_f$ years.

We define effectiveness as the proportion of cases averted
on the total cases possible under no intervention \cite{rodrigues2014cost}:
\begin{equation}
\label{eq:F}
\overline{F}=\frac{A}{t_f I(0)}=1-\frac{\int_0^{t_f}I^*(t)~dt}{t_f I(0)}.
\end{equation}
The total cost associated with the intervention  \cite{rodrigues2014cost} is
\begin{equation}
\label{eq:TCI}
TC=\int_0^{t_f} C \, \mathbbm{T}^*(t)I^*(t)~dt,
\end{equation}
where $C$ corresponds to the per person unit cost of the detection
and treatment of infectious individuals. Following
\cite{okosun2013optimal,rodrigues2014cost},
the average cost-effectiveness ratio is defined by
\begin{equation}
\label{eq:ACER}
ACER=\frac{TC}{A}.
\end{equation}
\begin{table}[!htb]
\centering
\caption{Sumary of cost-effectiveness measures. Parameters according 
to Table~\ref{tab:param} and $C=1.$}\label{tab:efficacy}
\begin{tabular}{c@{\hspace*{2.2cm}}c@{\hspace*{2.2cm}}c@{\hspace*{2.2cm}}c}
\toprule
$A$ \eqref{eq:A} & $TC$ \eqref{eq:TCI} & $ACER$ \eqref{eq:ACER} 
& $\overline{F}$ \eqref{eq:F}\\[1mm] \midrule
20.9 & 3459.1 & 165.5 & 0.00429\\ \bottomrule
\end{tabular}
\end{table}
Table \ref{tab:efficacy} summarizes the particular case we have analyzed.
These results evidence the limited effectiveness of the control variable
to diminish the HRSV infectious individuals.


\section{Conclusion}
\label{sec:conc}

Human Respiratory Syncytial Virus (HRSV)
is the most common cause of lower respiratory tract infection
in infants and children worldwide. In addition, HRSV causes serious
disease in elderly and immune compromised individuals. In this work,
we discussed a mathematical compartmental model for HRSV.
Estimation of parameters was done for real data
of Florida from September 2011 to July 2014,
minimizing the $l_2$ norm. The results show
that the proposed model fits well the reality under study.
Moreover, a sensitivity analysis was carried out for the basic
reproduction number, in the case when the transmission parameter
is taken to be the average value of $\beta(t)$ in the period under study,
proving that the most important parameters to have into account are
the rate of loss of infectiousness $\nu$ and
the average of transmission parameter $b_0$.
Our results from optimal control show that treatment has a limited
effect on HRSV infected individuals. This reinforces the importance
of developing a licensed vaccine for HRSV, which
is a subject under strong current development \cite{neuzil2016progress}.


\section*{Acknowledgements}

Rosa was supported by the Portuguese Foundation for Science and Technology (FCT)
through IT (project UID/EEA/50008/2013);
Torres by FCT through CIDMA (project UID/MAT/04106/2013)
and TOCCATA (project PTDC/EEI-AUT/2933/2014
funded by FEDER and COMPETE 2020). 

The authors are very grateful to the anonymous referees 
for a careful reading of the manuscript and for several 
questions and suggestions that improved the paper.



\end{document}